\documentstyle[12pt]{article}
\begin{document}

\title{ Continuum Cascade Model:  
  Branching Random Walk   for  Traveling Wave }
\author{ Yoshiaki Itoh \footnote {Institute of Statistical Mathematics,  
10-3 Midori-cho, Tachikawa, Tokyo,  190-8562  Japan, Email: 
 itoh@ism.ac.jp}  \\{\small The Institute of Statistical Mathematics and }\\{\small the Graduate University for 
Advanced Studies} 
\\{\small }
}

\date{\today}
\maketitle

 \begin{abstract}
  
 The cascade model generates random food webs. 
The continuum cascade model is a Poisson approximation of 
the cascade model. 
We have a simple nonlinear recursion for 
probability distribution of the longest 
chain length (the height) generated by the continuum cascade model.  
Assuming the 
traveling wave solution,  the   velocity selection principle for the Fisher-KPP 
equation works  for our 
recursion.  
 Here we have the recursion for the height of continuum cascade model 
  from the first passage time  of the
 left most  particle
of a branching Poisson point process. 
The  asymptotic probability distribution of 
 the height   is obtained   by a straightforward 
application of the Aidekon 
theorem for  the left most particle  of branching 
Poisson point process.  Hence    the traveling wave behavior 
is shown mathematically.


\end{abstract}



\section{Introduction}

 We introduced   a    nonlinear 
 recursion  \cite{ik} for the 
probability distribution of the longest chain length  (height of tree ) of 
the Poisson approximation of  a random  directed graph   model, cascade model of food web \cite{cbn, cn},  task graphs for parallel
processing \cite{gnpt} and Barak-Erdos graphs for stochastic 
order \cite{be, fk}.  
Assuming  a traveling 
wave solution, the  velocity selection principle  is naturally applied to our recursion  \cite{ik} as in the Kolmogorov-Petrovskii-Piskunov  argument for  the Fisher-KPP equation.  The asymptotic position of wave front of the constant velocity with the logarithmic correction term is obtained   \cite{ik}  by using an intuitive physical  argument \cite{km,mk} extending   the studies  on  the Fisher-KPP equation 
\cite{bramson,ds, ls, mckean, uchiyama}.  Here we obtain  the asymptotic probability distribution of the wave 
front mathematically.   
  The Aidekon  theorem \cite{aidekon, aidekon2, bdz} on branching random walk \cite{ar, hs, bk}, which is the analog of the Lalley and Selke theorem \cite{ls} of branching Brownian motion,  is  straightforwardly  applied to obtain the asymptotic probability distribution.  
Our recursion \cite{ik} gives  
the probability distribution of the minimum of the branching  Poisson point process.  
The solution to the Fisher-KPP  equation is given by   using 
a random shift, by derivative martingale of branching  Brownian motion,   of
the Gumbel distribution \cite{ls}.  The solution to our recursion is given by using 
a random shift, by derivative martingale of branching  Poisson point process,   of
the Gumbel distribution.

The cascade model \cite{cn}  generates a 
 food web at random.    
 Consider the random directed graph with vertex set $\{1, ... , n\}$ in which the 
${n\choose 2}$
directed edges $(i, j)$ with $i < j$ occur independently of each other with probability
$P = P_n = c/n$, and no edges with $i>  j$ occur. Such random graphs have been used
to model community food webs in ecology \cite{cn} and task graphs for parallel
processing in computer science \cite{gnpt}. The occurrence of an edge $(i, j)$  denotes, in
the biological context, that species $i$ is eaten by species $j$ or, in the computational
context, that task $i$ must be performed before task $j$. In both contexts, the
maximum path length is of interest.
Let us denote by $L = L_n$ the length (number of edges) of the longest (directed)
 path starting from vertex 1, and by $M = M_n$ the length of the longest path
(starting from any vertex). For a food web, $M$ represents the length of the longest
food chain. For a parallel computation in which each task takes one unit of time
and where the number of processors is sufficiently large, $M + 1$ represents the
processing time.

\begin{figure}[t]
\special{epsfile=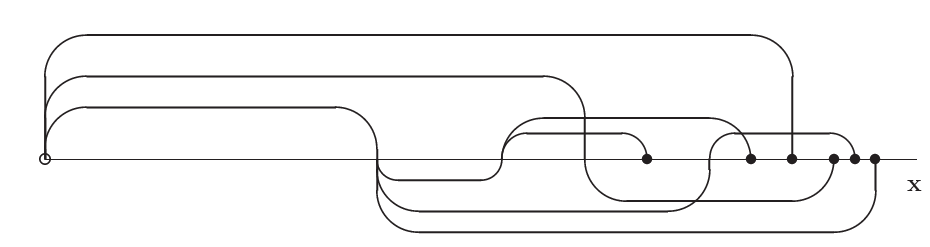 hscale=0.75  vscale=0.75} 
\vspace{3.5cm}
\caption{ Continuum cascade model,    \cite{ik}.} 
\label{continuum}
\end{figure}

  The probability  of  $k$ occurrences of 
  directed edges from the  vertex $1$ is given by the binomial distribution 
  ${n\choose k}p^{k}(1-p)^{n-1-k}$. 
   The longest chain length  starting from each end of the  occurred edges, which have the vertex 1 as the other end,   is not statistically   independent with  others   
and  we can not find the simple recursion to obtain the probability distribution of the longest chain length. 
However for  the Poisson approximation of the cascade model 
(continuum cascade model)  we have the statistical independence and have a very simple nonlinear recursion for the probability distribution 
of the longest chain length.   
Food webs typically involve a huge number of species, while the average predation per
species is usually not too large.  Hence, it is interesting to investigate large food webs with
$n\to \infty$, $c\to  0$ with finite  $np = x $, which gives the Poisson approximation of the orginal cascade model. 
 We are interested in the probability distribution  of the longest chain length $L_n$,  
starting
from vertex 1, as 
$n$ tends to infinity.  
Let us call  the Poisson approximation of the cascade model  the continuum cascade model \cite{ik}.
 In the illustrative picture
Fig. \ref{continuum} we draw only the vertex set and links from a  
tree,  generated by the continuum  cascade model,  initiating at the
origin (the open circle on the picture). Namely, we draw all links emanating from the
origin indicating direct predation on the basal species (there are 3 such predators in the
picture); then we draw all the links from these direct predators 
(4 such predators in Fig. \ref{continuum}); etc. Links are drawn in a cascade manner thereby explaining the name of the
model. 
Bottom preys are often called basal species. We reserve the term `basal species' only for the species at
the origin which, according to the definition of the continuum cascade model, can never be a predator
independently on the choice of links, \cite{ik}. 

\begin{figure}[t]
\special{epsfile=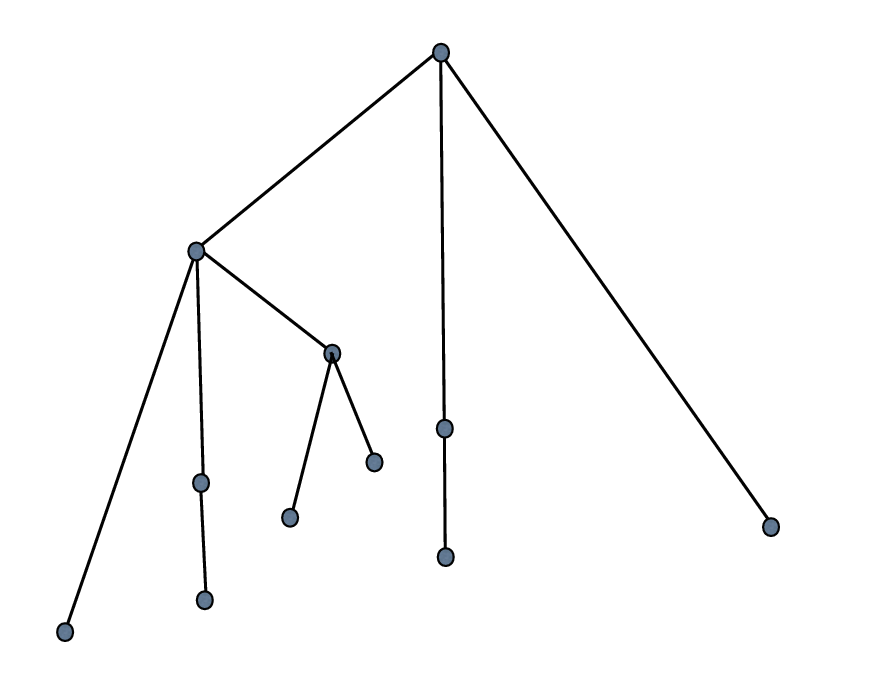 hscale=0.7  vscale=0.55} 
\vspace{6cm}
\caption{A tree generated by the continuum cascade model  with the basal species (the vertex at the top) playing the
role of the root. The height of this  tree is equal to 3, 
 \cite{ik}.} 
\label{tree}
\end{figure}
For the above illustrative picture Fig. \ref{continuum} we have  a tree with 10
links and 11 vertices. Six of these vertices (closed circles in Fig. \ref{continuum}) are terminal,
that is, there are no links emanating from them. 
It is convenient to utilize a more traditional way of plotting trees; the
cascade tree pictured Fig. \ref{continuum}  is presented on Fig. \ref{tree}. 
The
size,  the number of terminal vertices, the heght etc.  in  trees generated by 
the continuum cascade model
fluctuate from realization to
realization. In this paper we study the asymptotic  probability distribution of  the longest chain length (the height) of the continuum cascade model.

{\bf Remark 1}  The expected length of 
the longest chain length of  the original cascade model is obtained recently with an interesting correction term to the  constant $e$   \cite{fk, mr}.

{\bf Remark 2} It is pointed out \cite{gab} that 
 our continuum cascade model   is  also studied  under the name of 
Poisson weighted 
infinite 
 tree (PWIT), \cite{af, as}) in probabilistic combinatorial optimization.   
The closure of vertex 1 in the cascade model 
 converges in distribution to the PWIT  as $n$ tends to infinity \cite{gab}.

{\bf Remark 3}  In making  the continuum cascade model 
we extended  the ideas on 
random sequential bisection model  \cite{si,km,mk}  which is a continuum binary 
search tree. The random sequential bisection model gives 
an analogous asymptotic behavior  to the binary search tree  of $n$ keys for  the sorting algorithms  \cite{devroye, robson}. 
Applying the KPP velocity selection principle to  the Hattori and 
Ochiai  conjecture  \cite{ho1, ho2, km} for the random sequential bisection, 
the correction term  in the asymptotic expected height of the random sequential bisection is 
given   in \cite{km}.    
The correction term for  the expected height of binary search tree is   
 obtained   mathematically \cite{reed} ( see also \cite{drmota, szz}). The pick up stick model \cite{szz} analyzed by using the generating function makes a bridge between the random sequential  bisection model and the binary search tree.

\section{Nonlinear recursion for the height}

Now we consider the probability distribution of the the height,  which is 
the length (number of edges) of the longest (directed)
 path starting from vertex 1,
of the tree
 generated by our continuum cascade model. Let us define our model more precisely.

i) At step   1  we generate   $N_x$ points by  the Poisson distribution  
$Pr (N_x=k)=\frac{x^k}{k!}e^{-x}$ 
on the interval 
 $[0,x]$.  
Each of the  $N_x$ points is mutually independently distributed uniformly 
at random 
 on $[0,x]$. 

ii) At step  $u (>1)$,  for each generated point at $x-y$    
generated at the step  $u-1$,  generate $N_{x-y}$ points by the Poisson
 distribution $Pr (N_{x-y}=k)=\frac{(x-y)^k}{k!}e^{-(x-y)}$ 
uniformly at random on the interval $[y,x]$, 
  independently from  other points at step  $u$ and 
independently from the   points generated in the previous steps.  

iii) We make  $N_{x-y}$ edges from the point $y$ to each of the 
$N_{x-y}$ points generated on the interval  $[y,x]$.

iv)  We continue the above ii) and iii) recursively 
 as long as we have  at least one  new  point.  We stop the generation of points at step $H(x)$ when no new point is generated.

 The  $H(x)$ is the height of the tree 
   generated by the continuum cascade model  on 
 $[0,x]$. 
 Let
 \begin{eqnarray}
P_n (x)\equiv P(H(x)\leq n).
\end{eqnarray}
When  $k$ points, $x-y_1, x-y_2,...,x-y_k$ are generated at step 1, the probability,
 that 
the height is not larger than $n-1$,  is $P_{n-1} (y_1)P_{n-1} (y_2)\cdots P_{n-1}(y_k)$. 
 Since  each $y_i$ is distributed uniformly at random on $[0,x]$ and $k$ is 
 distributed by the Poisson  distribution,  
we have the following  recursion \cite{ik}   
for  the probability  $P_n (x)\equiv P(H(x)\leq n)$.   

For $n=0$, 
\begin{eqnarray}
P_n(x)=e^{-x}, \label{initial}
\end{eqnarray}
while  for $n\geq 1$, 
\begin{eqnarray}
\nonumber P_n(x)&=&  e^{-x}+ 
\sum _{k=1}^{\infty} \frac{x^k }{k!}e^{-x}  \frac{1}{x^k}  \int _0^x \cdots  
\int _0^x  P_{n-1}(y_1)\cdots P_{n-1}(y_k)\,dy_1\cdots dy_k
\\\nonumber
&=&  e^{-x}  \sum _{k=0}^{\infty} \frac{1}{k!} \left(   \int _0^x  P_{n-1}(y)\,dy \right)^k
\\&=&  \exp\left[-x + \int _0^x  P_{n-1}(y)\,dy)\right]. 
\label{eq:height2}
\end{eqnarray}

\section{Numerical traveling wave solution}
 We apply  the Aidekon theorem \cite{aidekon, aidekon2}
  to   show  the following observations  \cite{ik}  for 
 equation  (\ref{eq:height2})  mathematically 
 in later sections. 
 \begin{enumerate}
 \item 
 Numerical  traveling wave solution. 
Numerically, the probability distribution $P_n(x)$  has a traveling wave
shape with the width of the front remaining  finite as shown in Fig. \ref{tw}. 
\begin{figure}[t]
\special{epsfile=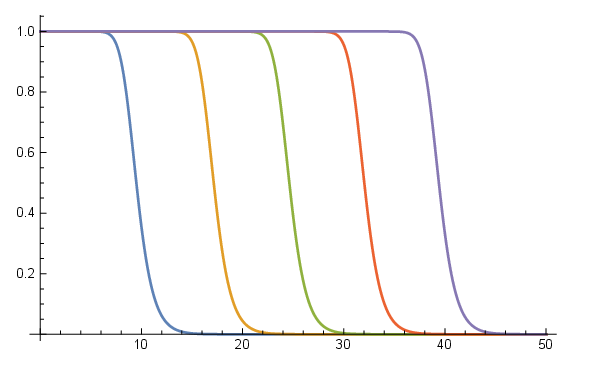 hscale=0.9  vscale=0.9} 
\vspace{6cm}
\caption{Traveling waves, 
the distribution $P_n(x)$ versus $x$ obtained by iterating the recursion. $P_n(x)$ is shown 
for $n = 20, 40, 60,  80, 100$  (left to right).  Iterations were performed  for 
a discrete approximation (\ref{discrete}) given in section \ref{numerical}. 
 } 
\label{tw}
\end{figure}
\item 
 Velocity selection. 
Assuming  an  approximation by a traveling wave form for  larger  
  $n\gg 1$, 
\begin{equation}
\label{TW}
P_n(x)\to \Pi(x-x_f), 
\end{equation}
with the  front position $x_f$ growing linearly with `velocity' equal to
 $e^{-1}$ \cite{ik}:
\begin{equation}
\label{front}
x_f\simeq vn, \quad v=\frac{1}{e}.
\end{equation}
 The velocity selection principle    \cite{kpp} 
gives $v = e^{-1}$ \cite{ik} by using an analogous argument  to the case of  
binary search tree \cite{km}.  We see that  the wave front $x_f$ should 
advance  asymptotically by a constant velocity 
$v = e^{-1}$,   from the  probabilistic argument for the 
 cascade model    \cite{cn, newman}.  
\item 
 Logarithmic correction.  An  analogy \cite{bramson, uchiyama} from the 
Fisher-KPP equation  gives a logarithmic correction to the front 
position  as
\begin{eqnarray}
x_f = \frac{n}{e} + \frac{3}{2e}\ln n + O(1). \label{front2}
\end{eqnarray}
It is convenient to think about $x$ and $n$ as space and time coordinates, 
so that the front of
the traveling wave was advancing. 

\item  Finite width of the front. 
The probability distribution $P_n(x)$  has asymptotically a traveling wave
shape with the width of the front remaining finite, and this is essentially equivalent to the finite
width of the height distribution.

\end{enumerate}
\section{Branching random walk for the solution of the recursion}\label{brw} 
We consider  the asymptotic behavior of   branching random
 walk  \cite{aidekon, aidekon2, ar,  bdz, bz, hs}.  We follow the notation and argument  by  Aidekon \cite{aidekon}. The process starts with one particle
 located at 0. At time 1, the particle dies and gives birth to a point
 process $\cal{L}$. Then, 
at each time $n \in N$, the particles of generation $n$ die and give birth 
to independent copies of
the point process $\cal{ L}$, translated to their position. If $\mathbf{T}$ is 
the genealogical tree of the process, 
we see that $\mathbf{T}$ is a Galton-Watson tree, and we denote by $|x|$ the 
generation of the vertex
$x \in \mathbf{T}$ (the ancestor is the only particle at generation 0).  For each 
$x \in \mathbf{T}$, we denote by
$V (x) \in  R$ its position on the real line. 
With this notation, ($V (x); |x| = 1$) is distributed as
$\cal{L}$.  The collection of positions ($V (x); x \in \mathbf{T}$) defines 
 our branching random walk.

We assume that we are in the boundary case \cite{hs}
\begin{eqnarray}
 E[\sum_{|x|=1}1]>1,~~\label{moment1}\\
E[\sum_{|x|=1}e^{-V(x)}]=1,~~\label{moment2}\\
E[\sum_{|x|=1}V(x)e^{-V(x)}]=0. \label{moment3}
\end{eqnarray}
 
 Every branching random walk satisfying mild assumptions can be reduced 
to this case by
some renormalization.  Notice that we may have
\begin{eqnarray}
\sum_{|x|=1} 1 = \infty \label{infty}
\end{eqnarray}
 with positive probability \cite{aidekon2}. We
are interested in the minimum at time $n$ 
\begin{eqnarray}
M_n := \min\{ V (x); |x| = n\},
\end{eqnarray}
where $\min |\emptyset|= \infty$. Writing for 
 $y \in  R \cup \{\pm \infty \}$,  $y_+ := \max(y,  0)$, we introduce 
the random variable
\begin{eqnarray}
X:=\sum_{|x|=1}e^{-V(x)},~~\\ \tilde{X}:=\sum_{|x|=1}V(x)_{+}e^{-V(x)}.
\end{eqnarray}
We assume  
 that
the distribution of $\cal{L}$ is non-lattice,
 we have
\begin{eqnarray} \label{moment4}
E[\sum_{|x|=1}V(x)^2e^{-V(x)}]<\infty
\end{eqnarray}
\begin{eqnarray}\label{moment5}
E[X(\ln_{+}X)^2]<\infty,~~
E[\tilde{X}(\ln_{+}\tilde{X})]<\infty
\end{eqnarray}
To state the  result, we
introduce the derivative martingale, defined for any $n>  0$ by
\begin{eqnarray}
 D_n :=\sum_{|x|=n}V (x)e^{-V (x)}.
\end{eqnarray}
From \cite{bk, aidekon2}
 ( Proposition A.3 in the  Appendix of \cite{aidekon2}), we know that the martingale 
converges
almost surely to some limit $D_{\infty}$, which is strictly positive on 
the set of non-extinction of $\mathbf{T}$.
Notice that under conditions (\ref{moment1}), (\ref{moment2}), (\ref{moment3}),  
the tree $\mathbf{T}$ has a positive probability 
to survive.

 There exists a constant $C^{*}>0$ such that for
 any real $x$,
\begin{eqnarray}\label{aidekon-theorem}
\lim_{n\to \infty}P(M_n>\frac{3\ln n}{2}+x)
=E[e^{-C^{*}e^x D_{\infty}}], 
\end{eqnarray}
(Theorem 1 in \cite{aidekon, aidekon2},  see \cite{bdz} for an
 elementary approach).

\section{Asymptotic probability for  the longest chain length}\label{minima}
The process starts with one particle
 located at 0. At time 1, the particle dies and gives birth to the point
 process $\cal{L}$, with intensity 1   on $[0,\infty)$. Then,
at each time $n \in N$, the particles of generation $n$ die and give birth 
to independent copies of
the point process $\cal{ L}$, translated to their position. 
At each time , we kill all particles to the right of $x$.   
Denote position of left-most particle
in this (extended) tree at $n$-th generation by $H_n$. Since
\begin{eqnarray}
P(H(x)\leq n-1)=P(H_n \geq x), 
\end{eqnarray}
we see 
\begin{eqnarray}\label{pnph}
P_{n-1}(x) = P(H_n \geq  x).
\end{eqnarray}
To   normalize  for equation (\ref{moment3})
we replace the original Poisson Point Process 
of intensity 1 on $[0,\infty)$  by the 
Poisson Point Process  of intensity $\frac{1}{e}$ on $[-1,\infty)$, as $\cal{L}$ in
 section \ref{brw}.
Then the  conditions (\ref{moment1}),  (\ref{moment2}),  (\ref{moment3})  hold.  
We see the inequality  (\ref{moment1}),  since  expected number of children here is infinite. 
For the  identity (\ref{moment2})  we  have 
\begin{eqnarray}
\int_{-1}^\infty e^{-y} \frac{dy}{e}  =1
\end{eqnarray}
and for the identity (\ref{moment3}) we have
\begin{eqnarray}
\int_{-1}^\infty y e^{-y} \frac{dy}{e} 
=0. 
\end{eqnarray}
The distribution of this Poisson point process  
$\cal{L}$ is non-lattice and
the  moment conditions  (\ref{moment4})  and  (\ref{moment5}) hold by
exponential decay (\ref{moment2}) for $V(x)$. We have 
\begin{eqnarray} 
\nonumber&& E[\sum_{|x|=1}V(x)^2e^{-V(x)}]
\\&=&\int_{-1}^{\infty} x^2 e^{-x}dx <\infty. 
\end{eqnarray}
The  total number of children is assumed to be 
finite almost surely in  \cite{aidekon}.  However the argument  \cite{aidekon} is  applied  to the above extension
to $[-1,\infty)$, as shown in \cite{aidekon2}.

The position $H_n$  of 
the  left-most particle at generation $n$ is 
given by  using  $M_n$ for 
\begin{eqnarray}
M_n = e H_n - n. \label{normalize}. 
\end{eqnarray}
Considering  equation (\ref{normalize}), 
\begin{eqnarray}
M_n>z+\frac{3}{2}\ln n, 
\end{eqnarray}
means  
\begin{eqnarray}
  H_n >\frac{z+n+\frac{3}{2}\ln n}{e}.
\end{eqnarray} 
Hence 
from equation (\ref{pnph}), 
\begin{eqnarray}
P (M_n>z+\frac{3\ln n}{2})&&=P(H_n>\frac{z+n+\frac{3}{2}\ln n}{e})\\
&&=P_{n-1}(\frac{z+n+\frac{3}{2}\ln n}{e}). 
\end{eqnarray} 
Put  $z/e=x$, then from 
   equation  (\ref{aidekon-theorem})  ( Aidekon  \cite{aidekon, aidekon2}),  for the solution $P_{n-1}$  to 
equation  ({\ref{eq:height2})
we have  
\begin{eqnarray}\label{distfront}
\lim_{n\to \infty} P_{n-1}(x+\frac{n}{e}+\frac{3}{2e}\ln n))
= E[\exp(-C^* e^{e x} D_\infty )], 
\end{eqnarray} 
 which gives the asymptotic probability  on   the  longest chain length
 (on the position of  wave front).

 For 
$x=0$ of equation (\ref{distfront}),  we have
\begin{eqnarray}\label{height4}
\lim_{n\to \infty}P_{n-1}(\frac{n}{e}+\frac{3}{2e}\ln n))=
                    E[\exp(-C^* D_\infty )], \label{frontapproach}
\end{eqnarray}
which should be less than 1 and larger than 0,  since $D_{\infty}$ is mathematically 
shown to be 
strictly positive \cite{aidekon, aidekon2, bk}.

\section{Numerical observation}\label{numerical}
\begin{figure}[t]
\special{epsfile=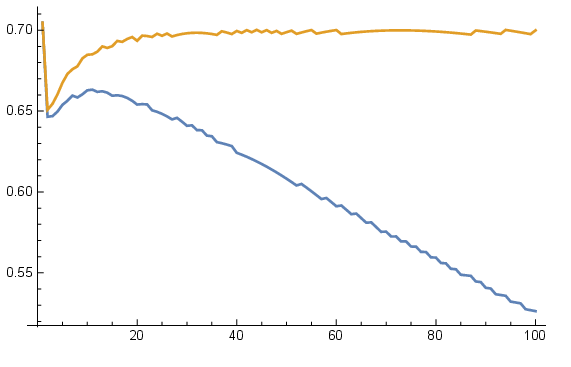 hscale=0.7  vscale=0.7} 
\vspace{5.cm}
\caption{$P_{n-1}(\frac{3}{2e}\ln n+\frac{n}{e}))$ for $n=1,2,...,100$
}\label{figfront}
\end{figure} 

Putting   $\tilde{x}\triangle$ for $x$, and  giving   the discrete initial value  
for   equation  (\ref{initial}), we 
consider  a recursion as a discretization of 
equation  (\ref{eq:height2})  
for $1\leq n$, 
\begin{eqnarray}\label{discrete}
f_n(\tilde{x})=\exp[-\tilde{x}\triangle +(\sum_{\tilde{y}=1}^{\tilde{x}} f_{n-1}(\tilde{y}))\triangle].
\end{eqnarray}
The 
numerical value  $f_n(\tilde{x})$    for $P_n(x)$  in 
Fig. \ref{tw}  and  Fig. \ref{figfront} are  obtained from  equation  (\ref{discrete})  for  
$\triangle=0.01$ by using the software  Mathematica.  The 
  numerical values  
\begin{eqnarray}
f_{n-1}(\frac{1}{\triangle}(\frac{3}{2e}\ln n+\frac{n}{e}))
\end{eqnarray}
for 
$P_{n-1}(\frac{3}{2e}\ln n+\frac{n}{e}))$, are shown  by the lower curve 
in Fig. \ref{figfront}.   
 Putting  
$e/\alpha$ instead of exponential $e$, the numerical values   
\begin{eqnarray}\label{alpha}
f_{n-1}(\frac{\alpha}{\triangle}(\frac{3}{2e}\ln n+\frac{n}{e}))
\end{eqnarray}
 for $\triangle=0.01$ and $\alpha=0.9855$ are
shown by the upper curve in 
Fig. \ref{figfront},  which  seems  to approach quickly to a constant.  We carried out 
  calculations and see  for example  
$\triangle=0.001$ and $\alpha=0.9977$ the  value of  (\ref{alpha})  quickly approaches
 to a constant.    
Our numerical calculations seem 
to suggest   $\alpha \to 1$ 
  as $\triangle \to 0$,  which supports  equation (\ref{frontapproach}) for 
the wavefront
 numerically.

{\bf Acknowledgements }  
The author thanks  Amir Dembo 
and Ofer Zeitouni for   suggesting  to  him   an essential idea given 
in   section \ref{minima}.  
The author thanks 
Elie Aidekon for    the 
 reply  \cite{aidekon2} to the question of Amir Dembo on \cite{aidekon},  
which is essential  for  the present paper.  
This work is   supported in part by US National Science Foundation 
Grant   DMS1225529 to Rockefeller University and JSPS Grant-in-aid for 
Scientific Research 
23540177.




\end{document}